\documentclass[11pt]{amsart}
\usepackage{amsmath}
\usepackage{amsthm}
\usepackage{amsbsy}
\usepackage{amsopn}
\usepackage{amsmath,amscd}
\usepackage{amssymb}
\usepackage{amsfonts}
\usepackage{multirow}
\usepackage{mathrsfs} 
\usepackage[official ]{eurosym}

\newcommand\cF{\mathcal{F}}

\def\bC{\mathbf{C}}

\def\bR{\mathbf{R}}

 \usepackage[all]{xy}
                \usepackage{amsthm,amsfonts,amsmath,amscd,amssymb,epsfig,verbatim,
}

\usepackage{graphicx}
\usepackage{epstopdf}
\date{\today} 
 
\usepackage{graphicx}

\usepackage[dvipsnames,svgnames,x11names,hyperref]{xcolor}
\usepackage{url,graphicx,verbatim,amssymb,enumerate,stmaryrd}
\usepackage[pagebackref,colorlinks,citecolor=Bittersweet,linkcolor=Bittersweet,urlcolor=Bittersweet,filecolor=Bittersweet]{hyperref}
\usepackage{tikz} 
\usetikzlibrary{arrows,decorations.pathmorphing,backgrounds,fit,positioning,shapes.symbols,chains,calc}
\usepackage[all]{xy}
\xyoption{matrix}
\xyoption{arrow}
\usepackage[hmargin=3cm,vmargin=3cm]{geometry}
\usepackage{setspace,kantlipsum}
\tikzset{help lines/.style={step=#1cm,very thin, color=gray},
help lines/.default=.5} 

\usepackage{booktabs}

\theoremstyle{plain}
\newtheorem{theorem}{Theorem}[section]

\newtheorem{lemma}[theorem]{Lemma}

\theoremstyle{definition}
\newtheorem{definition}[theorem]{Definition}
\newtheorem{inductive step}[theorem]{Inductive step}

\newtheorem{inductive lemma}[theorem]{Inductive Lemma}
\theoremstyle{remark}

\newtheorem{example}[theorem]{Example}

\newtheorem*{remark*}{Remark}

\newtheorem*{example*}{Example}

\newcommand{\wh}{\widehat}

\newcommand{\ve}{\varepsilon}

\numberwithin{figure}{section}

\title{Flexibility of singularities and beyond}

\author{Daniel  \'Alvarez-Gavela}
\address{Department of Mathematics \\ Brandeis University \\ Waltham, MA, 02453}
\email{dgavela@brandeis.edu}

\begin{document}

\maketitle
 \onehalfspacing
 
 \begin{abstract}
We survey a selection of Yasha Eliashberg's contributions to the philosophy of the h-principle, with a focus on the simplification of singularities and its applications.
 \end{abstract}
 
 \tableofcontents
 
 \section{Introduction}
 
 \subsection{Flexible mathematics}
 
{\em Flexible mathematics} is a loose term. Roughly speaking, a geometric problem is called {\em flexible} if in some sense it reduces to the underlying formal problem. The quintessential and foundational example is that of Hirsch-Smale-Whitney immersion theory. Here the geometric problem is to construct an immersion between manifolds $M \to N$, i.e. a smooth map $f:M \to N$ such that $df_x:T_xM \to T_{f(x)}N$ is an injective linear map for all $x \in M$. The underlying formal problem is that of constructing a bundle monomorphism $TM \to TN$, i.e. a bundle map $F: TM \to TN$ covering an {\em arbitrary} smooth map $f:M \to N$ such that $F_x:T_xM \to T_{f(x)}N$ is an injective linear map for all $x \in M$. When $\dim M < \dim N$ the Hirsch-Smale-Whitney theory states that the problem of constructing an immersion $M \to N$ is flexible: it is equivalent to the a priori much weaker problem of constructing a monomorphism $TM \to TN$ (in fact much more is true, as will be discussed in what follows).

If a problem is not flexible, then there is further geometry constraining the problem beyond what one would expect from the formal analogue, and we may call the problem {\em rigid}, or say that it exhibits {\em rigidity}. For example when $\dim M = \dim N$ and $M$ is closed (i.e. each connected component of $M$ is compact and boundaryless), existence of a bundle monomorphism $TM \to TN$ is a necessary but in general not sufficient condition for the existence of an immersion $M \to N$. So in this case the problem of existence of an immersion $M \to N$ is not flexible; it exhibits rigidity.

Under a narrative which Yasha has helped popularize, mathematics can be divided between these two categories: rigid and flexible. Some problems fall on one side of this dichotomy, and some problems fall on the other side. For example, most of algebraic geometry falls on the rigid side. It was a surprising discovery of 20th century mathematics that a good deal of geometric, symplectic and contact geometry falls on the flexible side, and moreover even some problems in complex and Riemannian geometry exhibit flexibility. Symplectic and contact geometry are distinguished in that they contain abundant phenomena on both sides of the rigid/flexible divide. Moreover, the boundary between rigidity and flexibility in symplectic and contact geometry is particularly fickle, with some problems walking a very fine line on the ridge that separates rigidity from flexibility (think existence of contact structures), and with some major open problems for which it is still unclear whether they belong to the rigid or to the flexible side.

A more recent narrative of Yasha's is that in fact every problem is flexible. Indeed, a possible objection to the above flexible/rigid dichotomy is that the formal problem underlying a given geometric problem often does not have an obvious formulation, and indeed is not a priori defined. One may sometimes move the goalposts by changing one's formulation of the underlying formal problem, perhaps far enough so that the original problem does indeed reduce to the modified formal problem, at which point we may declare the original problem to be flexible. From this viewpoint the task is to determine the optimal formulation for the formal problem. For example, one might think that the s-cobordism theorem is not a result in flexible mathematics because the problem of trivializing an h-cobordism does not reduce to the underlying homological problem, but it does reduce in some sense to the underlying K-theoretic problem. Indeed, if we take the vanishing of the Whitehead torsion as the formulation of the formal problem underlying the problem of trivializing an h-cobordism, then we may argue that the s-cobordism theorem is also an h-principle, after all.

The latest narrative of Yasha's is that there are two kinds of mathematicians, those who tell us what is not possible, and those who tell us what is possible. Roughly speaking, those who tell us what is not possible work on the rigid side: they create constraints and define invariants, thus putting an upper bound on the range of possible phenomena which may occur in a given geometric context. Those who tell us what is possible work in the opposite direction, constructing examples of phenomena which do occur in that given geometric context, or at least establishing their existence. In some sense the first group determines necessary conditions for the solution to a geometric problem to be possible, while the latter group establishes sufficient conditions. Often the two groups start quite far from each other, but after some progress they might inch closer and closer together, and if they're lucky they may eventually meet. 

Yasha Eliashberg has made substantial contributions to both the rigid and the flexible side of mathematics, sometimes telling us what is possible and sometimes telling us what isn't. This survey will focus on the flexible side. More specifically, this survey will focus on five specific frameworks of flexible mathematics developed by Eliashberg and illustrate them by telling the story of the simplification of singularities, from the Eliashberg viewpoint. The story will not be told exactly in chronological order, instead attempting to follow a self-contained mathematical narrative that builds up to the state of the art of present day research. The four theories are: (1) removal of singularities, (2) holonomic approximation, (3) surgery of singularities, (4) wrinkling, and (5) arborealization. There are some notable absences, including Yasha's work on overtwisted contact structures (covered in a separate essay in this volume by J. Etnyre) and Yasha's work on Stein and Weinstein manifolds (covered in a separate essay in this volume by K. Cieliebak). We will not dwell on these absences and instead refer the reader to the aforementioned essays.

One more word on terminology: it turns out that if a geometric problem does reduce in some sense to the underlying formal problem, then quite often a lot more is true. Notably, it is often the case that whenever existence of a formal solution (i.e. a solution to the formal problem) is sufficient for the existence of a genuine solution (i.e. a solution to the original problem), then one may in fact construct a genuine solution which is homotopic to any given formal solution, within the space of formal solutions. When this occurs, one says (following M. Gromov) that an h-principle holds (the ``h" is for {\em homotopy}). In fact in favorable circumstances it may even occur that the space of genuine solutions is (weakly) homotopy equivalent to the space of formal solutions, and it is common for even stronger statements to hold. The term {\em h-principle} is loosely used to describe problems which exhibit such flexible behavior, and one may use adjectives such as {\em parametric} or {\em relative} or {\em $C^0$-close} to further specify which flavor of the h-principle holds. 

\section{Removal of singularities}

While he was a PhD student, Yasha developed the technique of removal of singularities together with M. Gromov \cite{EG}. The technique of removal of singularities can be used to establish h-principles for a variety of classes of maps which in some sense avoid a given singularity type. The simplest class of non-singular maps consists of immersions, and indeed we will present the Hirsch-Smale-Whitney immersion theory as a showcase of the technique.

\subsection{Immersion theory}

The main result of the Hirsch-Smale-Whitney immersion theory \cite{Hi, Sm, W} is the following:

\begin{theorem}
If $\dim M < \dim N$ then the map $\text{imm}(M,N) \to \text{mon}(TM,TN)$ is a weak homotopy equivalence.
\end{theorem} 

Here $\text{imm}(M,N)$ is the space of immersions $M \to N$ and $\text{mon}(TM,TN)$ is the space of bundle monomorphisms $TM \to TN$. The map  $\text{imm}(M,N) \to \text{mon}(TM,TN)$ assigns to an immersion $f:M \to N$ its derivative $df:TM \to TN$, which is in particular a bundle monomorphism. 

Let us first focus on the zero-parametric existence result, namely that if $\dim M < \dim N$, then any monomorphism $TM \to TN$ is homotopic through such monomorphisms to the differential of an immersion $M \to N$. In particular, it follows that the existence of a monomorphism $TM \to TN$ is necessary and sufficient for the existence of an immersion $M \to N$. Put $m=\dim M<\dim N= n$.

\subsection{Main inductive step}

We consider first the case $N=\bR^n$. The data of a monomorphism $F:TM \to T\bR^n$ consists of a smooth map $f=(f_1,\ldots,f_n):M \to \bR^n$ (homotopically this is no data since $\bR^n$ is contractible) together with $n$ differential 1-forms $\varphi_1,\ldots,\varphi_n$ on $M$ such that at each point $x \in M$ the co-vectors $\varphi_1(x),\ldots, \varphi_n(x) \in T^*_xM$ span $T^*_xM$. This monomorphism $F$ is the differential of an immersion precisely when $\varphi_i = df_i$ for all $i=1,\ldots,n$. The idea is to modify $F$ one $\varphi_i$ at a time, so that we inductively have $\varphi_i = df_i$ for all $i \leq k$. At the end of the induction $k=n$ the construction of the immersion is complete, and checking that the underlying monomorphism is homotopic to the original monomorphism through monomorphisms is straightforward. Let us focus on the last stage of the induction process, which contains the main idea.

At the last stage of the inductive process the inductive hypothesis is that $\varphi_i = df_i$ for $i<n$, so at each $x \in M$ the co-vectors $df_1(x),\ldots, df_{n-1}(x), \varphi_n(x)$ span $T_x^*M$. Note therefore that the map $\widetilde{f}=(f_1,\ldots,f_{n-1}):M \to \bR^{n-1}$ always has differential of corank $\leq 1$, i.e. the kernel of $d \widetilde{f}$ is at most 1-dimensional. By the Thom-Boardmann theory of singularities \cite{B, AGV}, the locus $\Sigma^1=\{ x\in M : \ker(d \widetilde{f}_x ) \neq 0\} \subset M$ of such a map is generically a smooth codimension 1 submanifold of $M$. Note that along $\Sigma^1$ the kernel of $d\widetilde{f}$ is of dimension exactly 1, hence forms a line field $\ell$ along $\Sigma^1$. Thom-Boardmann theory also ensures that generically, the locus where $\ell$ is tangent to $\Sigma^1$ is a submanifold $\Sigma^{11} \subset \Sigma^1$ of codimension 1 in $\Sigma^1$ (hence codimension 2 in $M$). So $\ell \pitchfork \Sigma^1$ along $\Sigma^1 \setminus \Sigma^{11}$. More generally, Thom-Boardmann theory ensures the generic existence of a stratification $\Sigma^1 \supset \Sigma^{11} \supset \cdots \supset \Sigma^{1^{m-1}} \supset \Sigma^{1^m}$, where $1^k$ denotes a string of 1's of length $k$ and $\dim(\Sigma^{1^k}) = m-k$, such that $\Sigma^{1^{k+1}}$ is the locus of points in $\Sigma^{1^k}$ where $\ell$ is tangent to $\Sigma^{1^k}$, and hence $\ell \pitchfork \Sigma^{1^k}$ along $\Sigma^{1^k} \setminus \Sigma^{1^{k+1}}$. Finally, we note that the line field $\ell$ is trivial and indeed trivialized by $\varphi_n$. It is therefore not hard to inductively construct a function $f_n:M \to \bR$ using the stratification $\Sigma^1 \supset \Sigma^{11} \supset \cdots \supset \Sigma^{1^n}$ such that $df_n$ is nonzero on the line field $\ell$. Then $f=(f_1,\ldots,f_n):M \to \bR^n$ is the desired immersion.

\subsection{Conclusion of the argument} 

The previous stages of the induction are similar, with the caveat that a suitable non-integrable version of Thom-Boardmann theory must be used in order to ensure the existence of an appropriate stratification, see \cite{Ste} for a thorough discussion.  That being said, the argument goes through in the same way: starting with an arbitrary sequence of differential 1-forms $\varphi_1,\ldots,\varphi_n$ such that for all $x \in M$ the covectors $\varphi_i(x)$ span $T^*_xM$, we focus on the locus $\Sigma^1$ on which $(\varphi_2(x),\ldots,\varphi_n(x))$ do not span $T^*_xM$ and replace $\varphi_1$ with $df_1$ for a suitable function $f_1:M \to \bR$ by ensuring that $df_1$ is non-vanishing on the line field $\ell$ defined along $\Sigma^1$ on which $\varphi_2,\ldots,\varphi_n$ all vanish. Then we do the same for $(df_1,\varphi_3,\ldots,\varphi_n)$ and so on, until we get to the last stage as described above. 

With some book keeping one can convince oneself that the underlying monomorphism of tangent bundles is homotopic to the one we started with. Moreover, with some minor modifications one can arrange for the argument to work in the case of a general target manifold $N$ instead of $\bR^n$. Furthermore, the argument also works parametrically, as well as in relative form, and therefore gives a full proof of the main result of Hirsch-Smale-Whitney theory.

The technique of removal of singularities can be applied to a number of other settings, including to the problem of constructing embeddings \cite{Sz} (which a priori may not look like a problem of removing singularities).

\section{Holonomic approximation}

As the field of the h-principle matured, the explicit corrugation constructions that had been introduced in various forms in the work of Whitney, Hirsch, Smale \cite{Hi, Sm, W}, Nash \cite{Na} and others were successively abstracted into general frameworks, notably in the work of Gromov on flexible sheaves \cite{G} and convex integration \cite{G, Sp}. Following this vein, the holonomic approximation lemma was later formulated by Eliashberg and Mishachev \cite{EM1} as a versatile tool for proving h-principle type results in the presence of some positive codimension, in particular recovering Gromov's h-principle for open Diff-invariant partial differential relations on open manifolds. We will state the lemma and illustrate its use by giving another proof of the Hirsch-Smale-Whitney immersion theory.

\subsection{The holonomic approximation lemma}

To state it we recall $J^r(M,N)$ the space of $r$-jets of maps $M \to N$ (think Taylor polynomials of order $r$)  and the notion of a holonomic section $M \to J^r(M,N)$, which is a special type of section of the projection $J^r(M,N) \to M$ (basepoint of the Taylor polynomial) given by the Taylor polynomials of an actual map $f:M \to N$, i.e. $s(x)=j^rf(x)$, the order $r$ Taylor polynomial of $f$ at $x\in M$. For example, sections of $J^1(M,N)$ consist of bundle morphisms $F:TM \to TN$, i.e. an arbitrary smooth map $f:M \to N$ together with a family of arbitrary linear maps $F_x:T_xM \to T_{f(x)}N$. When $F_x=df_x$ the section is holonomic. 

\begin{figure}[h]
\includegraphics[height=3cm]{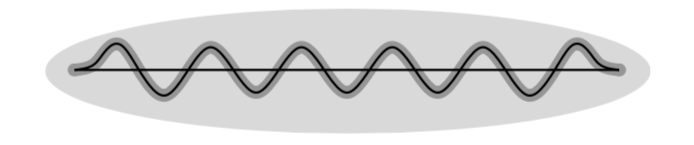}
\caption{A picture from the book {\em Introduction to the h-principle} \cite{EM} by Eliashberg and Mishachev showing a subset of positive codimension $K$, its wiggling $\varphi_1(K)$ and the neighborhood $U$ (darker shade) on which the holonomic approximation is defined.}
\label{fig: wiggling}
\end{figure}

Let $K \subset M$ be a reasonable subset of positive codimension (for example a subpolyhedron of a triangulation which contains no top dimensional simplices) and $s:M \to J^r(M,N)$ any section. The holonomic approximation lemma states the following (see Figure \ref{fig: wiggling} for an illustration of the wiggling).

\begin{lemma}[Holonomic approximation lemma] There exists an isotopy $\varphi_t:M \to M$, a neighborhood $U$ of $\varphi_1(K)$ and a holonomic section $\wh s : U \to J^r(U,N)$ such that $\wh s(x)$ is $C^0$-close to $s(x)$ for $x \in U$. \end{lemma}

For example, suppose that $r=1$ and $s:M \to J^1(M,N)$ is a section $s=(f,F)$ such that $F_x:T_xM \to T_{f(x)}N$ is injective for all $x \in M$. Then the lemma produces a holonomic section $\wh s = j^1g$ for a map $g:U \to N$ such that the linear map $dg_x:T_xM \to T_{g(x)}N$ is $C^0$-close to the linear map $F_x:T_xM \to T_{g(x)}N$ for all $x \in U$. If we take the approximation to be sufficiently $C^0$-close it follows that $dg_x$ is injective for all $x \in U$. If there exists an isotopy $\psi_t:M \to M$ such that $\psi_t(M) \subset U$ for $t\geq T$ we obtain an immersion $f=g \circ \psi_T:M \to N$. When $M$ is open, i.e. when no connected component is a closed boundaryless manifold, there always exists a polyhedron $K \subset M$ of positive codimension admitting such a $\psi_t$. Since furthermore the holonomic approximation lemma also holds in relative and parametric forms, this proves the following h-principle for immersions of open manifolds:

\begin{theorem}
If $M$ is open then the map $\text{imm}(M,N) \to \text{mon}(TM,TN)$ is a weak homotopy equivalence.
\end{theorem} 

Further, if $M$ is not necessarily open, but $\dim M < \dim N$, then to each injective bundle map $F_x:T_xM \to T_{f(x)}N$ we may assign its normal bundle $\nu_x=f^*(T_{f(x)}N/ F_x(T_xM) )$ and extend $F$ to an injective bundle map on the total space $V$ of the bundle $\nu \to M$. The manifold $V$ is open (of dimension $\dim V = \dim M)$ and it may be shrunk into an arbitrarily small neighborhood of the zero section of $\nu$, which is $M$. By applying the h-principle for open manifolds to $V$ as above and restricting to the zero section one obtains the full h-principle \cite{Hi,Sm,W} for immersions of (not necessarily open) manifolds into manifolds of greater dimension:

\begin{theorem}
If $\dim M < \dim N$ then the map $\text{imm}(M,N) \to \text{mon}(TM,TN)$ is a weak homotopy equivalence.
\end{theorem}

The above proofs touch on two big ideas in the philosophy of the h-principle, the first being Diff-invariance and the second being the micro extension trick. Indeed, the above proof of the h-principle for immersions of open manifolds can easily be adapted to prove the general Gromov h-principle for open Diff-invariant partial differential relations on open manifolds \cite{G,EM}. In particular the reader is invited to minimally adapt the above argument to prove the Phillips' h-principle for submersions of open manifolds \cite{Ph}, where we denote by $\text{sub}(M,N)$ the space of submersions $M \to N$ and by $\text{epi}(TM,TN)$ the space of bundle epimorphisms.

\begin{theorem}
If $M$ is open then the map $\text{sub}(M,N) \to \text{epi}(TM,TN)$ is a weak homotopy equivalence.
\end{theorem} 

The so-called micro extension trick used to deduce the h-principle for immersions of manifolds into manifolds of greater dimension from the h-principle for immersions of open manifolds cannot be used in the case of submersions, indeed no naive h-principle for submersions of closed manifolds holds, though we will return to submersions later once we allow wrinkles to enter into our life.

\subsection{Wiggling into the codimension}

We explain the proof of the holonomic approximation lemma for sections of the 1-jet bundle $J^1(\bR^2,\bR)$ defined on the cube $[-1,1]^2 \subset \bR^2$ with the subset $K$ equal to the $x$- axis $[-1,1] \times 0$. The data of a section is a family of linear polynomials $ax+by+c$ where $a,b,c$ are functions of $(x,y) \in [-1,1] \times [-1,1]$, but since we're working on an arbitrarily small neighborhood of the $x$-axis $K$ and may ignore a $C^0$-small error we will in fact ignore the $y$ coordinate and put $a(x)=a(x,0)$, $b(x)=b(x,0)$, $c(x)=c(x,0)$. 

Deform $K$ by a smooth isotopy supported in the unit cube such that $\varphi_1(K)$ is a sinusoidal curve of amplitude $\ve>0$ and period $4\delta$. On each quarter of the periodic oscillation, which happens on an interval $[k\delta, (k+1)\delta]$ we may define our holonomic section $\wh s$ on a neighborhood of $\varphi_1(K)$ by putting $\wh s = j^1g$ for $g(x,y)=a(\tau(y))x+b(\tau(y))y + c(\tau(y))$ where $\tau(y)$ is a function that increases from $k\delta$ to $(k+1)\delta$ as $y$ ranges from $-\ve$ to $\ve$. This serves as a piecewise definition for the function $g$ on some small neighborhood $U$ of $\varphi_1(K)$ and so we get a well-defined smooth function $g:U \to \bR$.

Let us compute $j^1g$. First, we consider the zeroth order part: $g(x,y)=a(\tau(y))x+b(\tau(y))y + c(\tau(y))$ is close to $a(x,y)x+b(x,y)y+c(x,y)$ because $(\tau(y),0)$ is close to $(x,y)$, indeed $\tau(y)$ is close to $x$ and $y$ is small. Hence the values of $a$, $b$ and $c$ at the points $(\tau(y),0)$ and $(x,y)$ are close to each other. Next, $\partial g/ \partial x = a(\tau(y) ) = a(a(\tau(y),0)$ which for the same reason as before is close to $a(x,y)$. Finally, note that $\partial g / \partial y $ has a term $b(\tau(y))$ which is again close to $b(x,y)$ for the same reason, plus the error term $$\tau'(y) \left(a'(\tau(y))x + b'(\tau(y))y + c'(\tau(y)) \right).$$ This error term can be ensured to be small since $\tau'$ is of order $\delta/\ve$ and the other term is bounded. Indeed, $\delta/\ve$ can be arranged to be arbitrarily small for example by taking $\delta=\ve^2$ and $\ve$ small. This completes the proof of the holonomic approximation lemma in this particular case.

The general case involves an inductive version of this argument, wiggling in each of the coordinate directions using the existing codimension, making the section holonomic one coordinate at a time. A simple calculation shows that to obtain the desired approximation in an $r$-jet bundle one needs to take $\delta, \ve>0$ small in the above construction so that $\delta/\ve^r$ is also arbitrarily small.

\section{Surgery of singularities}

\subsection{Maps between manifolds of the same dimension}

The technique of surgery of singularities constituted Yasha's PhD thesis. Let us give some context. At the time, Hirsch-Smale-Whitney immersion theory was established, and work by Poenaru and Phillips had led to the closely related h-principle for submersions of open manifolds.

However, when $M$ is closed and $\dim M=\dim N$ (so an immersion is the same as a submersion, and a monomorphism the same as an epimorphism), it is certainly false that the map $\text{imm}(M,N) \to \text{mon}(TM,TN)$ is a weak homotopy equivalence. Indeed we may even have $\text{imm}(M,N)$ empty while $\text{mon}(TM,TN)$ is nonempty, for example if we take $M$ to be a closed parallelizable manifold and $N$ to be Euclidean space. Hence when $\dim M= \dim N$ and $M$ is closed the problem of constructing immersions $M \to N$ is rigid; it falls outside of the realm of flexibility. What is then the best one may hope for in terms of flexibility? 

Well, the simplest kind of singularity for smooth maps between manifolds of the same dimension is the {\em fold}, which by definition is up to a change of coordinates given by the standard model $(x_1,x_2,\ldots,x_n) \mapsto (x_1^2,x_2,\ldots, x_n)$. Perhaps if we allow our maps to have folds we may optimistically hope to re-enter into the realm of flexibility.

If $S$ is a hypersurface of $M$ and $\dim M=\dim N$, we say that a map $f:M \to N$ is an $S$-immersion if $f$ is an immersion on the complement of $S$ and has fold type singularities along $S$. For simplicity we will assume in what follows that $S$ is two-sided. The formal analogue of an $S$-immersion $f:M \to N$ is a monomorphism $F:T_SM \to TN$, where $T_SM$, a vector bundle over $M$, denotes the result of cutting $TM$ open along $S$ and then re-gluing it along $S$ using an involution of the normal bundle of $S$. For example, if $M=S^n$ and $S=S^{n-1} \subset S^n$, the equator, then $T_{S^{n-1}}S^n$ is the trivial rank $n$ vector bundle over $S^n$.

Existence of a monomorphism $T_SM \to TN$ is a necessary condition for the existence of an $S$-immersion $M \to N$, indeed the differential $df:TM \to TN$ of an $S$-immersion $f:M \to N$ vanishes transversely along $S$, and its kernel is a line field which is transverse to $S$, so it is not hard to convince oneself that $df$ may be regularized to a monomorphism $\widetilde{df}:T_SM \to TN$, well defined up to contractible choice. Further, we have a space of $S$-immersions $\text{imm}_S(M,N)$, a space $\text{mon}(T_SM,TN)$ of bundle monomorphisms $T_SM \to TN$, and a well defined map up to homotopy $\text{imm}_S(M,N) \to \text{mon}(T_SM,TN)$, which assigns to an $S$-immersion $f:M \to N$ its (homotopically canonical) regularized differential $\widetilde{df}:T_SM \to TN$. What can one say about the map $\text{imm}_S(M,N) \to \text{mon}(T_SM,TN)$? Does an h-principle hold? 

Let us start with the earliest result in this direction. Suppose that there exists a bundle monomorphism $T_SM \to TN$, $M$ is connected, and $S$ is nonempty.  Let $T \subset M \setminus S$ be a two-sided hypersurface such that each connected component of $M \setminus S$ contains a connected component of $T$. An argument of Poenaru \cite{Po} shows that there exists a $V$-immersion $M \to N$, where $V$ consists of the disjoint union of $S$ and a finite union of $2k$ parallel copies of $T$ for some $k$. We briefly explain the argument.

The basic idea is to first use the h-principle for immersions of open manifolds on $M \setminus T$ to obtain an $S$-immersion $M \setminus T \to N$ (one should suitably prescribe the fold along $S$ first but the main point is that carving $T$ out gives us the codimension needed to make holonomic approximation or an equivalent tool work). Then once uses the 1-parametric h-principle for immersions of closed manifolds into manifolds of strictly greater dimension on a tubular neighborhood $U \simeq T \times (-1,1)$ of $T$, which we think of as a 1-parametric family of hypersurfaces $T_s=T \times \{s\}$, $-1\leq s \leq 1$. This produces a map $M \to N$ which is an $S$-immersion outside of $U$ and which restricts to an immersion on each parallel copy $T_s$ of $T$. 

Next, subdivide the parameter space $-1=s_0<\ldots <s_k=1$ so finely so that for each $s_j$ there exists a tubular neighborhood $U_j \to N$ of the immersion $T_{s_j} \to N$ such that the image of $T_{s_{j+1}} \to N$ is contained in the image of $U_j$ and moreover is graphical over $T_{s_j}$with respect to the collar coordinate of $U_j$. One may then further modify the 1-parametric family of immersions $T_s \to N$ by folding back and forth along these tubular neighborhoods, thus producing a smooth map $M \to N$ which has folds along the union $V$ of $S$ together with $2k$ parallel copies of $T$. This yields the desired $V$-immersion.

Whether the existence of a monomorphism $T_SM \to TN$ is sufficient for the existence of an $S$-immersion $M \to N$ (without additional folds) is not immediately apparent from Poenaru's folding argument, but can be deduced from the output of the Poenaru folding using surgery of singularities. Indeed, the technique of surgery of singularities allows for the simplification of the singularity locus in the source manifold under certain conditions. In particular we will explain a proof of the following existence h-principle for $S$-immersions, which is contained in Eliashberg's PhD thesis \cite{E1,E2}.

\begin{theorem}
If $M$ is connected and $S \subset M$ is nonempty, then map $\pi_0\text{imm}_S(M,N) \to \pi_0\text{mon}(T_SM,TN)$ is surjective.
\end{theorem}

\subsection{Direct and inverse surgeries}

Certain elementary surgeries of corank 1 singularities had appeared in the literature before the work of Eliashberg \cite{Le}. For example, consider the $\Sigma^{110}$ pleat, which for maps between surfaces is given by the normal form 
$$ (x,y) \mapsto (x,y^3+3xy)$$

\begin{figure}[h]
\includegraphics[height=5cm]{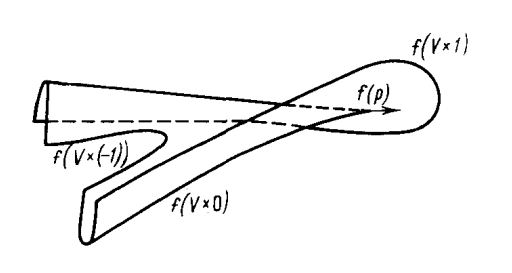}
\caption{A pleat with its characteristic vector, from Eliashberg's {\em On Singularities of folding type} 
\label{fig:characteristic-vector}
\cite{E1}.}
\end{figure}

The singular locus is the curve $\Sigma^1=\{x+y^2=0\}$, and the pleat occurs at the point $(0,0)$, with other points in the curve $x+y^2=0$ being folds (we recall that folds are denoted $\Sigma^{10}$ in the Boardman notation). The singular locus is contained in the half-space $x\leq 0$ bounded by the tangent line $T_{(0,0)}\Sigma^1$ to $\Sigma^1$ at the pleat point $(0,0)$. The characteristic vector $\nu$ at $(0,0)$ is uniquely determined up to contractible choice by demanding that it is nonzero and points into the other half-space $x \geq 0$, see Figure \ref{fig:characteristic-vector}. 

Suppose that $f:M \to N$ is a smooth map between surfaces with singular locus consisting of a curve $\Sigma \subset M$ and $p,q\in \Sigma$ are two pleat points. Suppose further that there exists an embedding $\alpha:[0,1] \to M$ with $\alpha(0)=p$, $\alpha(1)=q$, $\alpha(0,1) \cap \Sigma = \varnothing$, $\alpha \pitchfork \Sigma^1$ at $p$ and $q$, and with $\alpha'(0)$ and $-\alpha'(1)$ respectively the characteristic vectors $\nu_p$, $\nu_q$ for the pleats $p$ and $q$. Then it is possible to deform the map $f$ in a neighborhood $U$ of $\alpha([0,1])$ so that the two pleat points are cancelled against each other: more precisely the new singular locus $\widetilde \Sigma$ is equal to the union of $\Sigma \setminus U$ together with two arcs of $\Sigma^{10}$ folds which run parallel to $\alpha([0,1])$, and the open set $U$ no longer contains any pleats. This process results in a simultaneous Morse surgery of the $\Sigma^{110}$ and $\Sigma^{10}$ loci and is an example of direct surgery of singularities.

Similarly, it is straightforward to start with a smooth map between surfaces $f:M \to N$, $p \in M$ a $\Sigma^{10}$ fold point, and modify $f$ so that a pair of $\Sigma^{110}$ pleats are born along the curve $\Sigma^{10}$ near the point $p$. The characteristic vectors will point in opposite directions. This process also results in a simultaneous Morse surgery of the $\Sigma^{110}$ and $\Sigma^{10}$ loci and is also an example of direct surgery of singularities.

\subsection{Inverse surgeries via direct surgeries}

Direct surgeries are rather straightforward to realize. What is far less clear is how to achieve the reverse process, which is called inverse surgery of singularities. For example, suppose that $f:M \to N$ is a smooth map between surfaces, and there are two $\Sigma^{110}$ pleats on a $\Sigma^1$ curve which have only $\Sigma^{10}$ folds between them and whose characteristic vectors point in opposite directions. Can we deform $f$ so that the two pleats disappear, leaving only the curve of $\Sigma^{10}$ between them? This would be the reverse process to the direct surgery described in the last paragraph of the previous subsection. 

There difficulty is in some sense standard, indeed it is not always possible to cancel pairs of critical points of functions (for example for a function $S^1 \to \bR$), though creating cancelling pairs of critical points is always possible. Remarkably, it turns out that this rigidity only appears at the level of folds: all other inverse surgeries of corank 1 singularities can be geometrically realized. Eliashberg proved this by factoring (most) inverse surgeries as a composition of direct surgeries.

Let us illustrate the idea on the above concrete example. Again, we have $f:M \to N$ a smooth map between surfaces, and there are two $\Sigma^{110}$ pleats on a $\Sigma^1$ curve with an arc $A\subset \Sigma^1$ of $\Sigma^{10}$ folds between them and such that the characteristic vectors at the two pleat points in opposite directions. What one may in fact do in this situation is create (by a direct surgery) another pair of $\Sigma^{110}$ pleats, with characteristic vectors also pointing in opposite directions, on the same $\Sigma^{1}$ curve but just outside of $A$. One then cancels the four $\Sigma^{110}$ pleats against each other, using the other type of direct surgery described in the previous subsection. Note that in this cancellation one matches the two pairs of pleats with each other in a way that intermingles the two original pairs. This whole process all happens in a neighborhood $U$ of $A$ and the end result is a map $g:M \to N$, homotopic to $f$ rel. $M \setminus U$, which has $\Sigma^{10}$ folds on a curve $A'$ in $U$ which is isotopic to $A$ by an isotopy compactly supported in $U$. In particular the topology of the $\Sigma^1$ singular locus hasn't changed, though the two $\Sigma^{10}$ pleats have been removed. This concludes the process of geometrically realizing the inverse surgery of singularities in this particular case.

\subsection{Proof of the existence h-principle for $S$-immersions}

With most inverse surgeries factored as direct surgeries, Eliashberg obtained a number of results concerning the simplification of corank 1 singularities \cite{E1}. We concentrate on the above existence h-principle for $S$-immersions when $\dim M = \dim N$, though analogous results were also obtained in the general case $\dim M \geq \dim N$ \cite{E2}.

\begin{figure}[h]
\includegraphics[height=7cm]{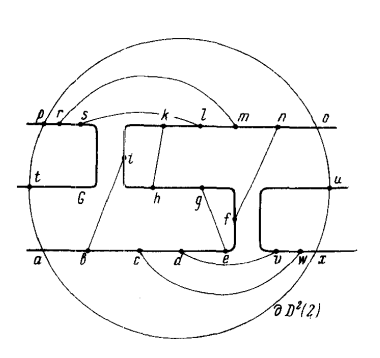}
\caption{A picture from Eliashberg's {\em On singularities of folding type} \cite{E1} indicating the sequence of surgeries needed to absorb a pair of concentric circles of folds into an existing locus. }
\label{fig: yashapic2}
\end{figure}

\begin{figure}[h]
\includegraphics[height=7cm]{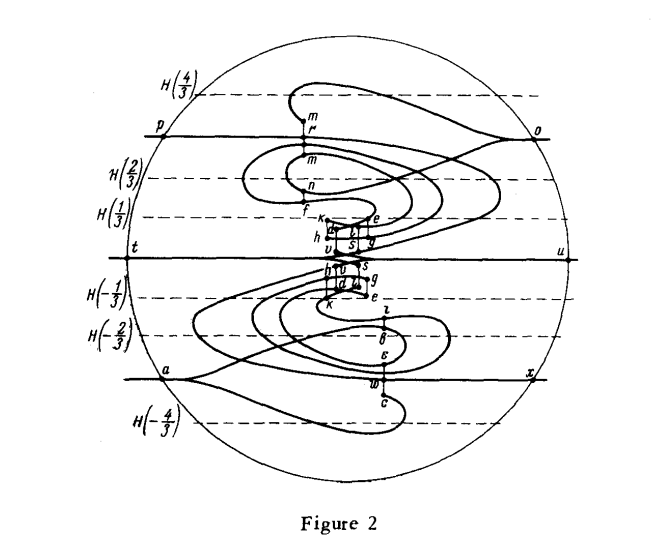}
\caption{Another picture from that same paper illustrating the effect of the surgeries in the target. }
\label{fig: yashapic1}
\end{figure}

We start with a monomorphism $T_SM \to TN$. For each component of $M \setminus S$ choose an embedding of $S^1 \times S^{n-2}$ into that component and call the disjoint union of all these hypersurfaces $T \subset M \setminus S$. By the Poenaru folding argument we may construct a $V$-immersion $M \to N$ for $V$ the disjoint union of $S$ and a bunch of parallel copies of $T$. One may then use surgery of singularities to absorb $T$ into $S$. Indeed, when $\dim M=2$ we may absorb two parallel circles of $\Sigma^{10}$ folds into a curve of $\Sigma^{10}$ folds by creating three pairs of $\Sigma^{110}$ pleats on each of the three curves using a direct surgery as described above and then cancelling them out against each other using the other type of direct surgery described above. The general case (i.e. arbitrary dimension) is similar, after multiplying everything by $S^{n-2}$, see Figures \ref{fig: yashapic2} and \ref{fig: yashapic1}. 

\section{Wrinkling}

\subsection{Wrinkling of mappings}

In the case of maps of smooth manifolds $f:M \to N$ with $\dim M = \dim N$, we have seen that allowing some kind of singularity is unavoidable. Of course the same is true when $\dim M \geq \dim N$. The theory of wrinkling \cite{EM2,EM3,EM4} shows that it is enough to allow one extremely simple type of corank 1 singularity locus, which is called a wrinkle, to obtain an h-principle type result, even in the absence of any positive codimension to help us. We will mostly restrict our discussion to the case $\dim M = \dim N$ for simplicity. In this setting a wrinkle of the map $f$ is a ball $B \subset M$ such that $f|_B$ is equivalent to the normal form
$$(x,y) \mapsto (x,y^3+3(\|x\|^2-1)y), \qquad (x,y) \in \bR^{n-1} \times \bR $$
Note that the normal form for the wrinkle has singularities on the unit sphere $S^{n-1} \subset \bR^n$, which consist of $\Sigma^{110}$ pleats on the equator $S^{n-2} \subset S^{n-1}$ and $\Sigma^{10}$ folds on the two hemispheres $S^{n-1} \setminus S^{n-2}$. 

For maps $f:M \to N$ with $\dim M > \dim N$ the model for the wrinkle is stabilized by a non-degenerate quadratic form $Q(z)$ of some index $0 \leq j \leq \left\lfloor (m-n)/2 \right\rfloor$.
$$(z,x,y) \mapsto (x,y^3+3(\|x\|^2-1)y+Q(z)), \qquad (x,y,z) \in \bR^{m-n} \times \bR^{n-1} \times \bR $$

\begin{definition} A map $f:M \to N$, $\dim M \geq \dim N$ is called a {\em wrinkled map}, if it is a submersion outside of a disjoint union of balls $B \subset M$ and each restriction $f|_B:B \to N$ is equivalent to the above normal form. \end{definition}

\begin{figure}[h]
\includegraphics[height=8cm]{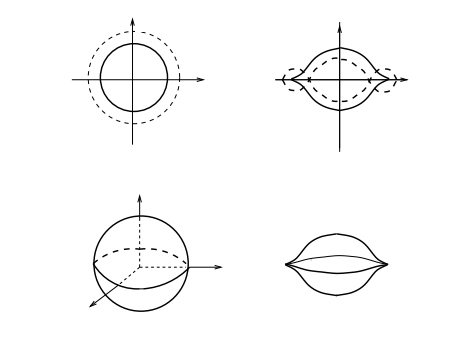}
\caption{Wrinkles in the source (left) and in the target (right), taken from the wrinkling saga \cite{EM2,EM3,EM4}. }
\label{fig: wrinkles}
\end{figure}

If $f:M \to N$ is a wrinkled map, then the differential $df$ of course degenerates along the wrinkles, but there is a homotopically canonical way to regularize it into a bundle epimorphism $\widetilde{df}:TM \to TN$. This is called the regularized differential of a wrinkled map. The existence form of wrinkled mappings theorem by Eliashberg and Mishachev \cite{EM2} states the following:
\begin{theorem} Any bundle epimorphism $F: TM \to TN$ is homotopic to the regularized differential of a wrinkled map $f:M \to N$. \end{theorem}

 The conclusion holds in relative form and also parametrically if one allows the wrinkles to be born and die in their standard embryo bifurcation. To be precise we may take any subset of the coordinates $x_j$ and consider them as parameters to obtain the normal form for fibered wrinkles. 
 
 \subsection{Soft and taut $S$-immersions}
 
 We now revisit the study of $S$-immersions $M \to N$ when $\dim M = \dim N$, $M$ is connected and $S$ is nonempty. It turns out that there is a dichotomy in that some $S$-immersions are flexible, while others are not. In fact the space $\text{imm}_S(M,N)$ decomposes as a disjoint union $\text{imm}_S(M,N) = \text{taut}_S(M,N) \coprod \text{soft}_S(M,N)$ where $\text{taut}_S(M,N)$, the space of {\em taut} $S$-immersions, consists of those $S$-immersions for which there exists an involution $M \to M$ which pointwise fixes $S$ and such the $S$-immersion is invariant under the involution. The space of {\em soft} $S$-immersions consists of those $S$-immersions which are not taut. 
 
Taut $S$-immersions can be thought of as rigid. For example, if $M \setminus S$ has two connected components $V_+$ and $V_-$, which are therefore exchanged by the involution, then $\text{taut}_S(M,N)$ is the product of the space of immersions $V_+ \to N$ (which abides by an $h$-principle since $V_+$ has nonempty boundary) and the space of diffeomorphisms of $V_+$ relative to the boundary $\partial V_+$ in a suitable sense.

Soft $S$-immersions are flexible in that an $h$-principle holds. Indeed, Eliashberg and Mishachev use the wrinkling technology to prove the following result when $\dim M \geq 2$ in \cite{EM5} (with the exception of the case where $\dim M=2$ and $N$ is closed, which to my knowledge is still open). 

\begin{theorem}
If $M$ is connected and $S$ is nonempty then the map $\text{soft}_S(M,N) \to \text{mon}_S(M,N)$ is a weak homotopy equivalence. 
\end{theorem}

Before we give the proof, let us first give a different proof of the surjectivity of the map $\pi_0 \text{imm}_S(M,N) \to \pi_0 \text{mon}(T_SM,TN)$, from a more wrinkled viewpoint. Given a bundle monomorphism $T_SM \to TN$, using the wrinkled mappings theorem it is not hard to construct a map $M \to N$ which has folds along $S$ and outside of $S$ is an immersion except for a finite disjoint union of wrinkles. One may then use surgery of singularities to absorb all these wrinkles into the existing fold locus $S$. At the last step one has two parallel $(n-2)$-spheres of $\Sigma^{110}$ pleats on the $(n-1)$-dimensional $\Sigma^1$ locus $S$ and one needs an inverse surgery to get rid of the pleats and end up with only $\Sigma^{110}$ folds along $S$. Fortunately this inverse surgery can be factored in terms of direct surgeries as explained in the discussion above. 

With this in mind, let's now try to understand what happens with this argument when parameters are introduced. Even with the addition of one parameter, one must encounter the embryo singularities which are the instances of birth/death of wrinkles. It is not immediately clear what to do with the above surgery argument at these bifurcation points. To overcome this difficulty it is convenient to use a slightly different approach, not surgering the wrinkles into $S$ but instead {\em engulfing} them into $S$. This is not possible in general, however it is possible for soft $S$-immersions. It turns out that soft $S$-immersions can be characterized by the presence of a local model, which is an instance of a notion that Eliashberg has popularized as a {\em virus of flexibility}. In the presence of a flexibility virus, the whole problem becomes globally flexible. In this case the virus is a zig-zag. 

\begin{figure}[h]
\includegraphics[height=4cm]{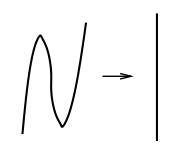}
\caption{A zig-zag, from \cite{EM5}.}
\label{fig: zig-zag}
\end{figure}

Given an $S$-immersion $M \to N$, a zig-zag is an embedding of a closed interval $A=[a,b]$ into $M$ which intersects $S$ transversely at two points and such that the composition of the embedding $A \to M$ and the $S$-immersion $M \to N$ is a smooth map $A \to N$ which factors as the embedding of an interval $B=[c,d]$ into $N$ and a smooth map between intervals $[a,b] \to [c,d]$ sending $a \mapsto c$ and $b \mapsto d$, which has exactly two non-degenerate critical points (a local maximum and a local minimum) in the interior of $[a,b]$, see Figure \ref{fig: zig-zag}. 

It is easy to see that if an $S$-immersion is taut, then it does not admit a zig-zag. Indeed the presence of the involution would force the arc $A$ to self-intersect. It is not too hard to convince oneself that the converse is also true: if an $S$-immersion does not admit zig-zags, then one may construct a suitable involution of $M$ by lifting paths from $N$ to $M$, and thus deduce that the $S$-immersion is taut. In conclusion: an $S$-immersion is soft if and only if it admits a zig-zag.

Once one has a zig-zag, one has as many zig-zags as one likes at one's disposal (take parallel disjoint arcs). One may then choose paths from the zig-zags to the location of embryonic birth/death of wrinkles and send the zig-zags along these paths to replace the wrinkles even before they are born, thus absorbing all the singularities coming from the wrinkling process into the fold locus $S$. This is the process known as engulfing. The key point is that soft $S$-immersions have sufficient flexibility to imitate the parametric wrinkling process using the already existing folds on $S$. 

The situation is reminiscent of contact structures, which also come in two types: {\em tight} and {\em overtwisted}. Overtwisted contact structures are characterized by the presence of a local flexibility virus and satisfy a full h-principle. Tight contact structures are those contact structures which are not overtwisted, and they are rigid. The development of the theory of overtwisted contact structures is certainly a highlight of Eliashberg's contribution to flexible mathematics, and is the subject of J. Etnyre's essay in this volume. 

\begin{figure}[h]
\includegraphics[height=4cm]{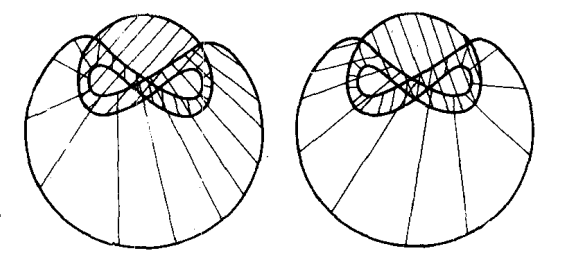}
\caption{The soft $S^1$-immersion $S^2 \to \bR^2$, illustrated as two immersions of the disk $D^2 \to \bR^2$ which agree on their boundary $\partial D^2 = S^1$, from \cite{E1}.}
\label{fig: soft}
\end{figure}

\begin{example}
There are two taut $S^1$-immersions $S^2 \to \bR^2$ up to homotopy through such, where $S^1 \subset S^2$ is the equator. These are: (1) the projection $S^2 \to \bR^2$, $(x,y,z) \mapsto (x,y)$ and (2) its post-composition with an orientation reversing diffeomorphism of $\bR^2$. According to the above theorem on $S$-immersions, there should also be only two soft $S^1$-immersions $S^2 \to \bR^2$ up to homotopy through such (and given one of them the other will be given by post-composition with an orientation reversing diffeomorphism of $\bR^2$). One may visualize a soft $S^1$-immersion $S^2 \to \bR^2$ as follows, in what is one of Yasha's favorite pictures. Start with the standard taut $S^1$-immersion $S^2 \to \bR^2$, i.e. the projection $(x,y,z) \mapsto (x,y)$. Use the direct surgery described above to create two cancelling pairs of cancelling $\Sigma^{110}$ pleats on the $\Sigma^1$ locus $S^1 \subset S^2$. Then use the other type of direct surgery described above to cancel the four $\Sigma^{110}$ pleats against each other, switching up the pairing as usual. Up to an isotopy of $S^2$ one obtains a new $S^1$-immersion $S^2 \to \bR^2$, which is soft. The projection of the fold locus to $\bR^2$ is illustrated in Figure \ref{fig: soft}. 
\end{example}

\subsection{Wrinkling of functions}

The parametric form of the wrinkled mappings theorem is very useful even when $N=\bR$, in which case we obtain applications to parametrized Morse theory. Indeed one can use the wrinkling technology to prove h-principles for functions $f:M \to \bR$ with mild singularities. By definition a function $f:M \to \bR$ with mild singularities is allowed to have Morse (quadratic) critical points $\sum_{i \leq k} x_i^2 - \sum_{i>k}x_i^2$ or Morse birth/death (cubic) critical points $x_1^3 + \sum_{1<i\leq k} x_i^2 - \sum_{i>k}x_i^2$, but nothing worse. 

The strongest result in this direction is formulated as follows. We define a {\em framed function} to be a function $f:M \to \bR$ with mild singularities which is decorated with the data of a framing of the negative eigenbundle of the Hessian of $f$ at each critical point of $f$. The space of framed functions is topologized so that the framings vary continuously in families and have to be suitably compatible at birth/death bifurcations. The main result from Eliashberg and Mishachev \cite{EM7}, following work of Igusa \cite{I1,I2} and Lurie \cite{Lu} is the following:

\begin{theorem} The space of framed functions on any manifold $M$ is contractible. \end{theorem}

The significance of the contractibility of the space framed functions is that it becomes possible to make homotopically canonical choices of (suitably decorated) Morse functions on smooth manifolds, which can be useful for geometric applications (see for example \cite{I3}).

\subsection{Wrinkling of embeddings}

The wrinkled mappings theorem is related but distinct to the wrinkled embeddings theorem. The wrinkled embeddings theorem is useful when one wants to simplify the tangencies of a submanifold $M \subset Y$ with respect to a foliation $\cF$ of $Y$, and is essential for many applications including those to parametrized Morse theory. If $\cF$ consists of the fibers of a fibration $\pi: Y \to N$, then tangencies of $M$ with respect to $\cF$ are the same as singularities of the map $\pi|_M: M \to N$, and our goal is to simplify these singularities. But of course we may only deform the smooth map $\pi|_M:M \to N$ through maps $M \to N$ of the form $\pi \circ \varphi_t$ for $\varphi_t:M \to Y$ an isotopy of $M$ in $Y$.

If $\dim M + \dim \cF < \dim Y$ (in the fibration case this means that $\dim M < \dim N$), then there is enough codimension available to deduce the h-principle from the holonomic approximation lemma or from any equivalent method of the h-principle arsenal.  However, when $\dim M + \dim F \geq \dim Y$, the most naive form of the h-principle certainly fails. For example even if $\gamma$ is homotopic to a distribution $\gamma'$ which is transverse to $M$, it will not be true in general that $M$ is isotopic to a submanifold $M'$ which is transverse to $\cF$. Some tangencies will be unavoidable. So the best one can do is to hope to simplify the tangencies as much as possible, both in terms of the model for the tangency as well as the topology of the tangency locus. And the simplest tangencies of them all are folds. 

An obvious necessary condition for $M$ to be isotopic to a submanifold $M' \subset Y$ such that the tangencies of $M'$ with respect to $\cF$ are all folds is that the distribution $\gamma = T \cF$ is homotopic to a distribution $\gamma'$ whose tangencies with respect to $M$ are all folds. The h-principle in this case says that this purely homotopy theoretic necessary condition is in fact also sufficient.

In order to prove such an h-principle in \cite{EM6}, Eliashberg and Mishachev pass through an intermediate object, called a wrinkled embedding. Let us for concreteness focus on the case $\dim Y = \dim M +1$. A wrinkled embedding $f:M \to Y$ is a topological embedding which is allowed to have singularities modeled on the normal form 
$$(x,y) \mapsto (x,y^3+3(\|x\|^2-1)y, \int_0^y(y^2+\|x\|^2-1)^2dz ), \qquad (x,y) \in \bR^{n-1} \times \bR$$

Note that the above formula consists of the normal form for a standard wrinkle together with an extra component which is a monotonically increasing in $y$ for all fixed $x$. Further, the partial derivative in the $y$ direction is
$$(x,y) \mapsto (0, 3(y^2+\|x\|^2-1), (y^2+\|x\|^2-1)^2) $$

\begin{figure}[h]
\includegraphics[height=4cm]{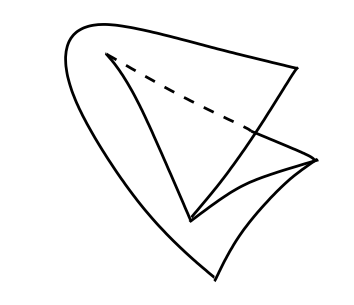}
\caption{One-half of an embedded wrinkle.}
\label{fig: half-wrinkle}
\end{figure}

So we observe that a wrinkled embedding has semi-cubic cusps above the fold points of the two hemispheres $S^{n-1} \setminus S^n$ which cancel in birth/death bifurcations along the equator $S^{n-2} \subset S^{n-1}$, see Figure \ref{fig: half-wrinkle}. 

Fibering along some subset of the $x=(x_1,\ldots, x_{n-1})$ coordinates we get the model for fibered wrinkled embeddings, in particular fibering over one coordinate $x_i$ gives the model for birth/death of wrinkles, with the model at the instant of bifurcation called an embryo. In particular in 1-parameter families of wrinkled embeddings $f_t:M \to Y$ we will allow wrinkles to be born and die along embryo bifurcations.

A wrinkled embedding $f:M \to Y$ has a well-defined Gauss map $G(df):M \to \text{Gr}_nY$, where $\text{Gr}_nY \to Y$ is the Grassmann bundle of $n$-planes in $Y$, which is given by $G(df)(x)=df_x(T_xM)$ at the smooth points and extended continuously to the singular locus by taking the limit. The additional flexibility provided by the wrinkles of wrinkled embeddings allows for the following remarkable statement: 

\begin{theorem} Given any homotopy $G_t:M \to \text{Gr}_nY$ of the Gauss map $G_0=G(dg)$ of a smooth embedding $g:M \to Y$, there exists a homotopy of wrinkled embeddings $g_t:M \to Y$ (i.e. we allow wrinkles to be born during the homotopy) such that $G(dg_t):M \to \text{Gr}_nY$ is arbitrarily $C^0$-close to $G_t$.
\end{theorem}

\begin{figure}[h]
\includegraphics[height=4cm]{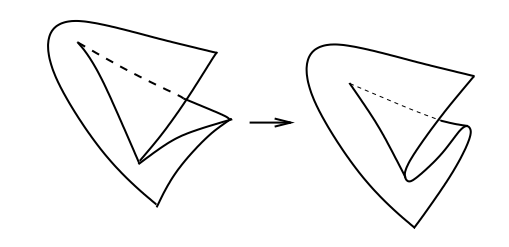}
\caption{Smoothing out (regularization) of (one-half of) an embedded wrinkle.}
\label{fig: regularization}
\end{figure}

One could use holonomic approximation to prove such a statement without wrinkles when $\dim M + \dim \cF < \dim Y$, but when $\dim M + \dim \cF \geq \dim Y$ the wrinkles are essential. Further, there is a homotopically canonical way of smoothing out the wrinkles of a wrinkled embedding, see Figure \ref{fig: regularization}, which essentially consists of replacing the normal form 
$$(x,y) \mapsto (x,y^3+3(\|x\|^2-1)y, \int_0^y(y^2+\|x\|^2-1)^2dz )$$
by the smooth embedding
$$(x,y) \mapsto (x,y^3+3(\|x\|^2-1)y, y ).$$

If a wrinkled embedding is transverse to a foliation $\cF$, then with a bit of care it is possible to perform the above regularization so that the smoothed embedding only has fold type tangencies with respect to $\cF$. This allows one to deduce an h-principle for the simplification of tangencies of smooth embeddings with respect to ambient foliations.

\subsection{Further applications}

The wrinkled embeddings theorem has enjoyed applications beyond the problem of simplifying tangencies, in particular to symplectic and contact geometry. Before briefly discussing these applications, we recall that the question of whether a section $s:M \to J^1(M,\bR)$ is holonomic is closely related to whether a submanifold $\Lambda \subset J^1(M,\bR)$ is isotropic, which by definition means that the 1-form $\alpha = \lambda - dz$ vanishes on $\Lambda$, where $\lambda$ is the canonical Liouville 1-form on $T^*M$ and we write $J^1(M,\bR) = T^*M \times \bR$ with $z$ the $\bR$ coordinate. Indeed if $\Lambda$ is graphical over $M$ then the two conditions are equivalent. 

When $\dim \Lambda < \dim M$, there is some codimension available so holonomic approximation (or an equivalent tool) can be used to prove an h-principle for isotropic embeddings. However when $\dim \Lambda = \dim M$, in which case an isotropic submanifold $\Lambda$ is called Legendrian, there is no codimension available and the most naive form of the h-principle fails. 

However, the wrinkled embeddings theorem can be suitably applied in the front space $J^0(M,\bR)=M \times \bR$ to prove an existence h-principle for Legendrian embeddings. Parameters pose difficulties, but in Murphy's PhD thesis \cite{Mu}, which was written under the supervision of Eliashberg, a full h-principle was proved for the class of loose Legendrians, which just like soft $S$-immersions are characterized by a flexibility virus. And just like for soft $S$-immersions, this flexibility virus allows for the engulfing of the parametric wrinkling process. These ideas are closely related to the development of flexibility in the theory of Weinstein manifolds \cite{CE}, another major contribution of Eliashberg to flexible mathematics, which is the subject of an essay by K. Cieliebak in this volume.

\subsection{Universal holes}

A wrinkle can be thought of as a way to fill in a {\em universal hole}. This is achieved by introducing the simplest possible singularities. One can in this way salvage a number of h-principles which do not hold in general but do hold for geometric structures which are allowed to have wrinkles. Whether or not anything can be salvaged without the introduction of wrinkles is a subtler question. For example, the h-principle for overtwisted contact structures \cite{E3, CPP, BEM} can be thought of as a construction that can fill a universal hole in contact geometry without introducing any singularities.

\section{The arborealization program}

We conclude this survey with a discussion of a flexible program which constitutes current work in progress by Eliashberg in collaboration with D. Nadler and the author. First, we will discuss the theory of singularities of wavefronts and their simplification. Then we will discuss the arborealization program for skeleta of Weinstein manifolds, which was initiated by Nadler \cite{N1,N2} and has also seen contributions from Starkston \cite{Sta}. Finally we will explain the relation between the two. 

\subsection{Singularities of wavefronts}

We briefly recall some foundational definitions from the theory of singularities of caustics and wavefronts, which like much of modern symplectic and contact geometry originates in the work of Arnold and his collaborators \cite{Ar}. 

A {\em symplectic manifold} is an even dimension smooth manifold $X$ equipped with a non-degenerate closed 2-form $\omega$, called a {\em symplectic form}. Non-degeneracy means that $\omega^n$ is non-vanishing, where $\dim X = 2n$. A {\em contact manifold} is an odd dimensional smooth manifold $Y$ equipped with a maximally non-integrable hyperplane field $\xi \subset TY$. If we assume for simplicity that that $\xi$ is co-orientable, so that $\xi = \ker(\alpha)$ for some 1-form $\alpha$, then the maximal non-integrability of $\xi$ amounts to the condition that $\alpha \land (d \alpha)^n$ is non-vanishing,  where $\dim Y = 2n+1$. Such a 1-form $\alpha$ is called a {\em contact form}.

A smooth embedding $f:L^n \to (X^{2n},\omega)$ (resp. $f:L^n \to (Y^{2n+1},\alpha)$) is called {\em Lagrangian} (resp. Legendrian) if $f^*\omega=0$ (resp. $f^*\alpha = 0$ in the co-orientable case, or $df(TL) \subset \xi$ in general). The image of a Lagrangian (resp. Legendrian) embedding is called a Lagrangian (resp. Legendrian) submanifold. 

A {\em Lagrangian fibration} is a fiber bundle $\pi : X \to B$ such that each fiber $\pi^{-1}(b) \subset X$, $b \in B$ is a Lagrangian submanifold. It is a standard theorem that every Lagrangian fibration is locally equivalent to a cotangent bundle $T^*B \to B$, where the symplectic form on $T^*B$ is given by $\omega = d\lambda$ for $\lambda$ the canonical Liouville 1-form. 

A {\em Lagrangian map} is a map $g:L^n \to B^n$ between manifolds of the same dimension such that $g=\pi \circ f$ for $f:L^n \to (X^{2n},\omega)$ a Lagrangian embedding and $\pi:X^{2n} \to B^n$ a Lagrangian fibration. Locally it suffices to understand the case $X=T^*B$, in which case every Lagrangian embedding $f:L \to T^*B$ can be locally lifted to a Legendrian embedding $\wh f: L \to J^1(B,\bR) = T^*B \times \bR$, where the contact form on $T^*B \times \bR$ is given by $dz-\lambda$ for $z$ the coordinate on the $\bR$ factor. 

It is useful to consider the {\em front projection} of such a lift, which is by definition the composition of $\wh f$ with the forgetful map $p:J^1(B,\bR) \to J^0(B,\bR) = B \times \bR$, i.e. $p$ is the product of the cotangent bundle projection $T^*B \to B$ and the identity $\bR \to \bR$. The image of the front projection $p \circ \wh f :L \to B \times \bR$ is sometimes called the {\em wavefront}. One may generically recover a Lagrangian embedding from its wavefront. 

Another standard model for the front projection is the map $S^*B \to B$ where $S^*B$ is the cosphere bundle of $B$. We recall that the cosphere bundle $S^*B$ is the quotient $T^*B/ \bR^+$ where $\bR^+$ acts by positive dilation of covectors, i.e. the positively projectivized cotangent bundle. The cosphere bundle $S^*B$ is equipped with a canonical contact structure and the composition of a Legendrian map $f:L^{n-1} \to S^*B$ with the cosphere bundle projection $S^*B \to B$ is also called a wavefront (note $\dim S^*B = 2n-1$ whereas $\dim J^1(B,\bR) = 2n+1$, where $\dim B = n$).

One may consider the problem of simplifying the singularities of a Lagrangian map with the additional constraint of only allowing deformations through Lagrangian maps. Equivalently, one can consider the problem of simplifying the singularities of wavefronts within the class of wavefronts. As in the case of smooth maps, the generic singularities of Lagrangian maps are impossible to understand, however the problem is flexible in that the strongest h-principle for the simplification of singularities of wavefronts that one could hope for does in fact hold. 

Concretely, M. Entov developed the theory of surgery of singularities in the category of Lagrangian maps \cite{EN}, and the author developed the theory of wrinkled emebddings in the category of Lagrangian maps \cite{AG2}, which also involved establishing some quantitative refinements of the holonomic approximation lemma \cite{AG1}. Both results constituted PhD theses under the supervision of Eliashberg. Thus, the problem of simplifying the singularities of wavefronts is fully flexible.  

\subsection{The arborealization program}

Weinstein manifolds are a distinguished class of open symplectic manifolds which can be thought of as the symplectic underpinning of Stein manifolds, i.e. holomorphic submanifolds of complex affine space $\bC^N$. A {\em Weinstein domain} is a compact manifold with nonempty boundary that can completed to a Weinstein manifold (the converse is true for Weinstein manifolds of finite type). 

A Weinstein domain $W$ is equipped with the following structure: an exact symplectic form $\omega = d \lambda$ with the choice of primitive $\lambda$ such that the vector field $Z$ which is $\omega$-dual to $\lambda$ is outwards pointing along $\partial W$ (a convexity condition), and a function $\phi:W \to \bR$ that has $\partial W$ as a regular level set and for which $Z$ is gradient-like (a taming condition for $Z$). One must be somewhat careful with the meaning of gradient like when $\phi$ is not Morse but we will not dwell in the details in this survey, see \cite{C}. Weinstein manifolds and Weinstein domains are important objects of study in symplectic topology, for example one encounters them often in mirror symmetry. 

A Weinstein domain (or manifold) always has the homotopy type of a CW complex of dimension no greater than $n= \frac{1}{2} \dim W$. Indeed, when $\phi$ is Morse it is not hard to verify that the index of a critical point of $\phi$ can be no greater than $n$, and then one can take the aforementioned CW complex to consist of the union $K$ of the stable manifolds of the critical points of $\phi$, which is a stratified subset with strata of dimension $\leq n$. In fact, these stable manifolds are isotropic.

In general the {\em skeleton} of a Weinstein domain $W$ consists of the subset $K= \bigcap_{t>0}Z^{-t}(W)$, where $Z^{-t}:W \to W$ denotes the flow of the vector field $-Z$ for time $t>0$. Up to homotopy of Weinstein structures, $W$ is completely determined by an arbitrarily small neighborhood of $K$. However, $W$ is in no reasonable sense determined by the stratified subset $K$, as can be seen even in the simple examples where $\phi$ has only two critical points, consisting of a minimum and an index $n$ critical point. The fundamental issue is that in general $K$ is highly singular, indeed too singular for $W$ to be recovered from $K$. 

However, there is a class of Lagrangian singularities, introduced by Nadler, which has the remarkable property that for a skeleton $K$ of a Weinstein manifold $W$ with singularities in this class, the Weinstein structure of $W$ is indeed determined up to deformation by the skeleton $K$, together with the discrete data of an {\em orientation structure} \cite{AGEN1}. Furthermore, these singularities are characterized locally up to contractible choice of symplectomorphism by combinatorial data \cite{AGEN1}. This should be thought of as an analogue of how open Riemann surfaces of finite type are determined up to deformation by a finite trivalent graph equipped with a cyclic ordering of the half-edges incident at each vertex. This distinguished class of Lagrangian singularities is called the class of {\em arboreal singularities}, due to the fact that they admit a natural indexing by (discretely decorated) rooted trees. 

\begin{figure}[h]
\includegraphics[height=4.3cm]{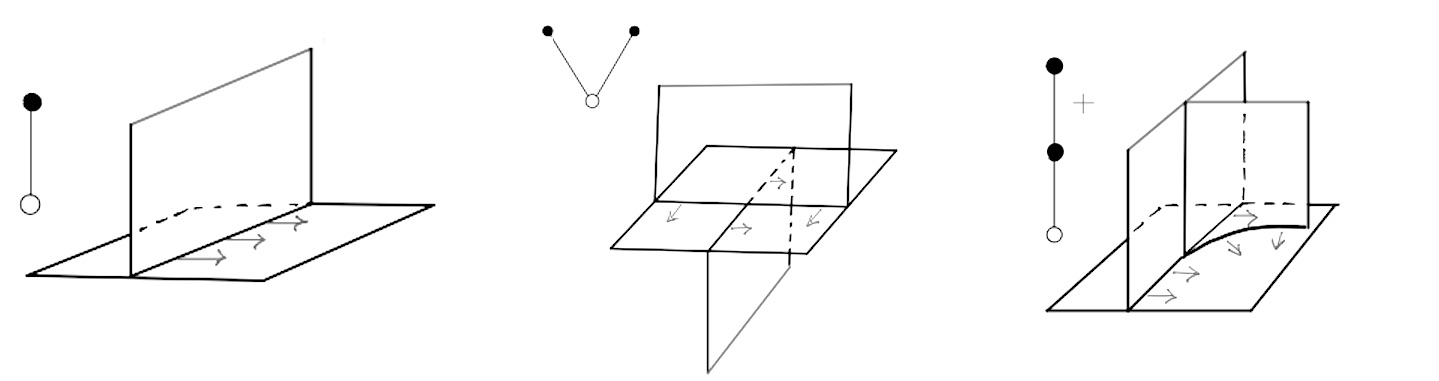}
\caption{Some arboreal singularities.}
\label{fig: arboreal}
\end{figure}

The arborealization program aims to reduce the study of the symplectic topology of Weinstein manifolds up to deformation to the study of the differential topology of arboreal spaces up to some standard Reidemeister type moves. By an arboreal space we mean a topological space which is locally modeled on arboreal singularities and is equipped with an orientation structure. However, there exist homotopy theoretic obstructions to a Weinstein manifold admitting an arboreal skeleton in its Weinstein deformation class. This condition is quite close to the existence of a polarization, i.e. a global field of Lagrangian planes, or equivalently the condition that the tangent bundle $TW$ is isomorphic to a bundle of the form $E \otimes \bC$ for $E$ a rank $n$ real vector bundle on $W$ (the symplectic structure on $W$ makes $TW$ a complex vector bundle). In fact, the existence of a polarization is necessary for the existence of a skeleton with singularities in the subclass of {\em positive} arboreal singularities. It turns out that this condition is also sufficient, which is an h-principle type result:

\begin{theorem} A Weinstein manifold $W$ admits a polarization if and only if it admits a positive arboreal skeleton. \end{theorem}

This result was obtained in joint work of Eliashberg with Nadler and the author \cite{AGEN3}. Current work in progress aims to establish a uniqueness counterpart to this existence result: namely that any two positive arboreal skeleta corresponding to the same polarized Weinstein manifold can be related by a finite set of Reidemeister type moves. 

\subsection{Relation between arborealization and flexibility of caustics} 

Let us focus on the case in which $W$ is a Weinstein domain with $\phi: W \to \bR$ a Morse function having only two critical points, a minimum $x_0$ and an index $n$ critical point $x_n$. The stable manifold $U$ of $x_n$ is an $n$-disk which intersects the boundary $\partial B$ of a standard Darboux ball $B$ about $x_0$ along a Legendrian $(n-1)$-sphere $\Lambda \subset \partial B$. The skeleton $K$ of $W$ in this case consists of the smooth $n$-disk $V=U \setminus B$ together with the radial cone of $\Lambda$ in $B$ centered at $x_0$, which can be highly singular at the point $x_0$. 

A first step towards spreading this singularity out would be to deform the Weinstein structure on $B$, in which $Z$ is the outwards pointing radial vector field, and replace it with the Weinstein structure on the cotangent bundle of a disk $T^*D^n$, in which $Z$ is the fiberwise outwards pointing radial vector field (so $Z$ now has zeros all long the zero section $D^n$ instead of just at the center of the ball $B$, and $\phi$ becomes Morse-Bott with $D^n$ as a critical submanifold). This can be achieved in such a way that $\Lambda$ becomes a Legendrian link in the cosphere bundle $S^*D^n$ of $D^n$, and the new skeleton of $W$ consists of the same smooth $n$-disk $V$ as before, together with the {\em fiberwise} radial cone of $\Lambda$ in $T^*D^n$, and together with the zero section $D^n$ of $T^*D^n$. Note that the new skeleton is singular along the image of the front $\Lambda \to D^n$, and is a smooth Lagrangian submanifold elsewhere.

When $n$ is small, the generic singularities of the front $\Lambda \to D^n$ are not so bad, and so one has successfully simplified the singularities of the skeleton (this was Starkston's approach to arborealization in the case $n=2$). However when $n>2$ the generic singularities of wavefronts are more complicated, and get arbitrarily bad as $n$ becomes larger and larger. Fortunately, the problem of simplifying the singularities of wavefronts abides by an h-principle as discussed earlier in this section, and so one may hope to simplify the singularities of the front $\Lambda \to B$. One may verify that existence of a polarization gives you precisely the homotopy theoretic condition needed to apply the h-principle, thus enabling the arborealization of the skeleton of $W$.

The general case in which $\phi$ has many critical points presents several additional difficulties. First, it is it not clear how a global homotopy theoretic hypothesis may be used to ensure the applicability of the h-principle for the simplification of singularities of wavefronts at each stage. To deal with this issue, the h-principle for the simplification of singularities of wavefronts must be strengthened to what we called an h-principle {\em without pre-conditions} in which we are always able to simplify the singularities of wavefronts at the expense of introducing certain combinatorial singularities, called ridges \cite{AGEN2}, which can then be easily arborealized. Second, one must moreover control the interaction of three or more strata, since after using the above scheme to appropriately fix the interaction of two strata one no longer has freedom to fix the interaction with a third or other strata. For this purpose the notion of positivity ended up playing a key role. 

The uniqueness theorem for positive arboreal skeleta up to Reidemeister moves presents even further difficulties, but we are hopeful that a satisfactory result will be attained, and on a good day are optimistic that interesting and useful applications to the symplectic topology of Weinstein manifolds will ensue. In the meantime, Yasha keeps having fun with his collaborators drawing pretty pictures, as he always has.

 \end{document}